\documentclass[12pt]{amsart}
\usepackage{euler, amsfonts, amssymb, latexsym, epsfig,epic}

\usepackage{amscd,amssymb}
\usepackage{verbatim}
\usepackage{pstricks}
\usepackage{pst-node}

\input{diagrams}

\usepackage{bibentry}

\usepackage{amsmath}
\usepackage{amsthm}
\usepackage{amssymb}
\usepackage{amsfonts}
\usepackage{amsxtra}
\usepackage{amscd}
\usepackage{epsfig}
\usepackage{verbatim}

\usepackage{latexsym,amstext,epsfig}
\usepackage[all, knot]{xy}
\xyoption{arc}

\newsymbol\pp 1275

\newcommand{\Hom}{\operatorname{Hom}}
\newcommand{\End}{\operatorname{End}}
\newcommand{\Ext}{\operatorname{Ext}}
\newcommand{\ext}{\operatorname{ext}}
\newcommand{\rep}{\operatorname{rep}}
\newcommand{\Proj}{\operatorname{Proj}}
\newcommand{\SI}{\operatorname{SI}}
\newcommand{\SL}{\operatorname{SL}}
\newcommand{\GL}{\operatorname{GL}}
\newcommand{\PGL}{\operatorname{PGL}}
\newcommand{\ZZ}{\mathbb Z}

\newcommand{\RR}{\mathbb R}
\newcommand{\NN}{\mathbb N}
\newcommand{\QQ}{\mathbb Q}
\newcommand{\coker}{\operatorname{coker}}

\newcommand{\Id}{\operatorname{Id}}

\newcommand{\supp}{\operatorname{supp}}
\newcommand{\dd}{{\operatorname{\mathbf{dim}}}}

\newcommand{\cone}{\operatorname{Cone}}
\newcommand{\rel}{\operatorname{relint}}

\newcommand{\filt}{\operatorname{filt}}
\newcommand{\twist}{\operatorname{twist}}

\newtheorem{theorem}{Theorem}[section]
\newtheorem{proposition}[theorem]{Proposition}

\newtheorem{lemma}[theorem]{Lemma}

\theoremstyle{definition}
\newtheorem{definition}[theorem]{Definition}
\newtheorem{remark}[theorem]{Remark}

\newtheorem{example}[theorem]{Example}

\newcount\cols
{\catcode`,=\active\catcode`|=\active
\gdef\Young(#1){\hbox{$\vcenter
{\mathcode`,="8000\mathcode`|="8000
\def,{\global\advance\cols by 1 &}%
\def|{\cr
      \multispan{\the\cols}\hrulefill\cr
       &\global\cols=2 }%
  \offinterlineskip\everycr{}\tabskip=0pt
  \dimen0=\ht\strutbox \advance\dimen0 by \dp\strutbox
    \halign
    {\vrule height \ht\strutbox depth \dp\strutbox##
      &&\hbox to \dimen0{\hss$##$\hss}\vrule\cr
     \noalign{\hrule}&\global\cols=2 #1\crcr
     \multispan{\the\cols}\hrulefill\cr%
   }
}$}} }

\begin{document}
\pagestyle{plain}

\title{Cluster fans, stability conditions, and domains of semi-invariants}
\author{Calin Chindris}
\address{University of Iowa, Department of Mathematics, Iowa City, IA 52242, USA}
\email[Calin Chindris]{calin-chindris@uiowa.edu}

\date{November 9, 2008; Revised: \today}

\bibliographystyle{plain}
\begin{abstract}
We show that the cone of finite stability conditions of a quiver $Q$ without oriented cycles has a fan covering given by (the dual of) the cluster fan of $Q$. Along the way, we give new proofs of Schofield's results \cite{S2} on perpendicular categories. From our results, we recover Igusa-Orr-Todorov-Weyman's theorem from \cite{IOTW} on cluster complexes and domains of semi-invariants for Dynkin quivers. For arbitrary quivers, we also give a description of the domains of semi-invariants labeled by real Schur roots in terms of quiver exceptional sets.
\end{abstract}

\maketitle

\section{Introduction} Given a quiver $Q$ without oriented cycles, the set of \emph{almost positive real Schur roots} of $Q$ is the set
$$
\Psi(Q)_{\geq -1}=\{\beta \mid \beta \text{~is a real Schur root~} \} \cup \{-\gamma_i \mid i \in Q_0\},
$$
where $\gamma_i$ is the dimension vector of the projective indecomposable representation at vertex $i$. For example, when $Q$ is a Dynkin quiver, the real Schur roots are precisely the positive roots of the corresponding Dynkin diagram. However, in general the set of real Schur roots has a much more complicated structure.

The (possibly infinite) cluster fan $C(Q)$ on the ground set $\Psi(Q)_{\geq -1}$ consists of the rational convex polyhedral cones generated by the compatible subsets of $\Psi(Q)_{\geq -1}$.  The details of our notations can be found in Section \ref{reco.quiver} and Section \ref{cluster-fan-stability.sec}.

Our goal in this paper is to give an interpretation of $C(Q)$ in terms of the geometry of the representations of $Q$. Following Ingalls-Thomas \cite{IT}, the cone $\mathcal{S}(Q)$ of \emph{finite stability conditions} is, by definition, the set of all $\sigma \in \QQ^{Q_0}$ for which there are finitely many $\sigma$-stable representations up to isomorphism.

Let $I: \QQ^{Q_0} \to \QQ^{Q_0}$ be the isomorphism defined by $I(\alpha)=\langle \alpha, \cdot \rangle_{Q}$ where $\langle \cdot, \cdot \rangle_Q$ is the Euler form of the quiver $Q$. Now, we can state our first result:

\begin{theorem} \label{stabilitycone-thm} Let $Q$ be a quiver without oriented cycles. Then $\mathcal{S}(Q)$ has a fan covering given by $\{I(\cone(C)) \mid C \text{~is a compatible subset of~} \Psi(Q)_{\geq -1} \}$.
\end{theorem}

To prove the theorem above we use techniques from quiver invariant theory, developed mainly by Derksen and Weyman \cite{DW1, DW2}, King \cite{K} and Schofield \cite{S2}. Using the $\sigma$-stable decomposition for dimension vectors and the $A_{\infty}$-formalism, we give a new proof of Schofield's Embedding Theorem \cite{S2} which plays a fundamental role in our study:

\begin{theorem}(see also \cite[Theorem 2.5]{S2}) \label{S-embed-thm} Let $\alpha$ be a pre-homogeneous dimension vector and let $\sigma$ be either $\langle \alpha, \cdot \rangle$ or $-\langle \cdot, \alpha \rangle$.
\begin{enumerate}
\renewcommand{\theenumi}{\arabic{enumi}}
\item There are finitely many, up to isomorphism, $\sigma$-stable representations $E_1, \dots, E_l$ with $l \leq |Q_0|-1$. Moreover, the $E_i$ are exceptional representations.

\item If $\beta_i=\dd{E_i}$ then after rearranging $\mathcal{E}=(\beta_1, \dots, \beta_l)$ is a quiver exceptional sequence.

\item Let $Q(\mathcal{E})$ be the quiver with vertices $1, \dots, l$, and $-\langle \beta_i, \beta_j \rangle$ arrows from $i$ to $j$ for all $1 \leq i \neq j \leq l$. Then there exists an equivalence of categories from $\rep(Q(\mathcal{E}))$ to $\rep(Q)^{ss}_{\sigma}$ sending the simple representation $S_i$ of $Q(\mathcal{E})$ at $i$ to $E_i$. Consequently, if
    $$I:\NN^{Q(\mathcal{E})_0}=\NN^l \rTo \NN^{Q_0}$$ is defined by
$$
I(\eta(1),\dots,\eta(l))=\sum_{i=1}^l\eta(i)\beta_i,
$$
then

$$
\langle \eta, \gamma \rangle_{Q(\mathcal{E})}=\langle I(\eta), I(\gamma)\rangle_{Q}.
$$
\end{enumerate}
\end{theorem}

In \cite{IOTW}, Igusa et al. initiated the study of cluster fans via domains of semi-invariants of quivers. In fact, their motivation was two-fold since domains of semi-invariants are also related to the Igusa-Orr \cite{IO} pictures from the homology of nilpotent groups. Let us briefly recall the definition of domains of semi-invariants (for further details, see Section \ref{reco.quiver}). If $\beta$ is a dimension vector of $Q$, the domain of semi-invariants $D(\beta)$ is defined by
$$
D(\beta)=\{\alpha \in \QQ^{Q_0} \mid \langle \alpha, \beta \rangle =0 \text{~and~} \langle \alpha, \beta' \rangle \leq 0, \forall \beta' \hookrightarrow \beta \}.
$$

It was proved in \cite[Theorem 8.1.7]{IOTW} that for a Dynkin quiver $Q$, the $(|Q_0|-1)$-skeleton of its cluster fan can be covered by the domains of semi-invariants labeled by the real Schur roots of $Q$. This result can also be obtained directly from our Theorem \ref{stabilitycone-thm}. In fact, we can show:

\begin{theorem}\label{domains-cluster-thm} Let $Q$ be a connected quiver without oriented cycles. Then $Q$ is either a Dynkin quiver or a generalized Kronecker quiver if and only if
$$
\bigcup_{\beta}D(\beta)= \bigcup_{C} \cone(C),
$$
where the union on the left is over all real Schur roots $\beta$ while the union on the right is over all compatible sets $C$ with (at most) $|Q_0|-1$ elements.
\end{theorem}

In order to describe the domains of semi-invariants for arbitrary quivers, we need to work with quiver exceptional sets instead of compatible ones.

\begin{theorem}\label{domains-exceptional-thm} Let $Q$ be a quiver without oriented cycles and let $\beta$ be a real Schur root. Then
there are finitely many quiver exceptional sets $\mathcal{E}_1, \dots, \mathcal{E}_m$, each of size at most $|Q_0|-1$, such that
$$
D(\beta)=\bigcup_{1 \leq i \leq m} \cone(\mathcal{E}_i).
$$

Consequently,
$$
\bigcup_{\beta} D(\beta)=\bigcup_{\mathcal{E}} \cone(\mathcal{E}),
$$
where the union on the left is over all real Schur roots $\beta$ while the union on the right is over all quiver exceptional sets $\mathcal{E}$ of cardinality at most $|Q_0|-1$.
\end{theorem}

The layout of this paper is as follows. In Section \ref{reco.quiver}, we recall the main tools from quiver invariant theory. This includes King's criterion for semi-stability of quiver representations, Derksen-Weyman's First Fundamental Theorem, the Saturation Theorem, and the Reciprocity Property for semi-invariants of quivers. Schofield's results on perpendicular categories are reviewed in Section \ref{Embed.sec} where we give new proofs of his results (see Theorem \ref{S-embed-thm} and Theorem \ref{S-preh-semi-inv-thm}). Cluster fans and stability conditions for quivers are discussed in Section \ref{cluster-fan-stability.sec} where we also prove Theorem \ref{stabilitycone-thm} and Theorem \ref{domains-cluster-thm}. In Section \ref{domains-exceptional.sec}, we study domains of semi-invariants via exceptional sets and prove Theorem \ref{domains-exceptional-thm}.

\section{Recollection on quiver invariant theory}\label{reco.quiver} In this section, we review the main tools from quiver invariant theory that will be used in the latter sections. Let $Q=(Q_0,Q_1,t,h)$ be a finite quiver with vertex set $Q_0$ and arrow set $Q_1$. The two functions $t,h:Q_1 \to Q_0$ assign to each arrow $a \in Q_1$ its tail \emph{ta} and head \emph{ha}, respectively.

Throughout this paper, we work over an algebraically closed field $K$ of characteristic zero. A representation $V$ of $Q$ over $K$
is a collection $(V(i), V(a))_{i \in Q_0, a\in Q_1}$ of finite-dimensional $K$-vector spaces $V(i)$, $i \in Q_0$, and $K$-linear maps $V(a)\in \Hom_{K}(V(ta), V(ha))$, $a \in Q_1$. The dimension vector of a representation $V$ of $Q$ is the function $\dd{V}: Q_0 \to \ZZ$ defined by $(\dd{V})(i)=\dim_{K} V(i)$ for $i\in Q_0$. A dimension vector $\alpha \in \ZZ^{Q_0}_{\geq 0}$ is said to be \emph{sincere} if $\alpha(i)>0$ for all $i \in Q_0$. Let $S_i$ be the one-dimensional simple representation at vertex $i \in Q_0$ and let us denote its dimension vector by  $\varepsilon_i$.

Given two representations $V$ and $W$ of $Q$, we define a morphism $\varphi:V \rightarrow W$ to be a collection of $K$-linear maps
$(\varphi(i))_{i \in Q_0}$ with $\varphi(i)\in \Hom_K(V(i), W(i))$, $i \in Q_0$, and such that $\varphi(ha)V(a)=W(a)\varphi(ta)$ for all $a \in Q_1$. We denote by $\Hom_Q(V,W)$ the $K$-vector space of all morphisms from $V$ to $W$. Let $V$ and $W$ be two representations of $Q$. We say that $V$ is a subrepresentation of $W$ if $V(i)$ is a subspace of $W(i)$ for all $i \in Q_0$ and $V(a)$ is the restriction of $W(a)$ to $V(ta)$ for all $a \in Q_1$. In this way, we obtain the abelian category $\rep(Q)$ of all quiver representations of $Q$.

\textbf{From now on, we assume that our quivers are without oriented cycles.} Let $P_i$ be the projective indecomposable representation at vertex $i \in Q_0$ and let us denote its dimension vector by $\gamma_i$; we call $\gamma_i$ a \emph{projective root}.

A representation $V$ is said to be a \emph{Schur} representation if $\End_{Q}(V) \cong K$. We say that $V$ is a \emph{rigid} representation if $\Ext^1_Q(V,V)=0$. Finally, we say that $V$ is an \emph{exceptional} representation if $V$ is a rigid Schur representation. The dimension vector of a Schur representation is called a \emph{Schur root} while the dimension vector an exceptional
representation is called a \emph{real Schur root}. For example, the projective roots are real Schur roots.

Given two representations $V$ and $W$ of $Q$, we have the Ringel's \cite{R} canonical exact sequence:
\begin{equation}\label{can-exact-seq}
0 \rightarrow \Hom_Q(V,W) \rightarrow \bigoplus_{i\in
Q_0}\Hom_{K}(V(i),W(i)){\buildrel
d^V_W\over\longrightarrow}\bigoplus_{a\in
Q_1}\Hom_{K}(V(ta),W(ha)),
\end{equation}
where $d^V_W((\varphi(i))_{i\in Q_0})=(\varphi(ha)V(a)-W(a)\varphi(ta))_{a\in Q_1}$ and $\Ext^1_{Q}(V,W)=\coker(d^V_W)$.

The Euler form of $Q$ is the $\ZZ$-bilinear form on $\ZZ^{Q_0}$ defined by
\begin{equation}\label{Euler-prod}
\langle \alpha,\beta \rangle_{Q} = \sum_{i\in Q_0}
\alpha(i)\beta(i)-\sum_{a \in Q_1} \alpha(ta)\beta(ha).
\end{equation}
Of course, we can extend this bilinear form to $\RR^{Q_0}$. (When no confusion arises, we drop the subscript $Q$.)

It follows from (\ref{can-exact-seq}) and (\ref{Euler-prod}) that
\begin{equation} \label{Euler-formula}
\langle\dd{V}, \dd{W} \rangle=
\dim_K\Hom_{Q}(V,W)-\dim_K\Ext^1_{Q}(V,W).
\end{equation}

\subsection{Semi-invariants and semi-stable representations} For a given dimension vector $\beta$ of $Q$, the representation space
of $\beta-$dimensional representations of $Q$ is defined by $$\rep(Q,\beta)=\bigoplus_{a\in Q_1}\Hom_{K}(K^{\beta(ta)}, K^{\beta(ha)}).$$
If $\GL(\beta)=\prod_{i\in Q_0}\GL_{\beta(i)}(K)$ then $\GL(\beta)$ acts algebraically on $\rep(Q,\beta)$ by simultaneous conjugation, i.e., for $g=(g(i))_{i\in Q_0}\in \GL(\beta)$ and $W=(W(a))_{a
\in Q_1} \in \rep(Q,\beta)$, we define $g \cdot W$ by $$(g\cdot W)(a)=g(ha)W(a)g(ta)^{-1}, \forall a \in Q_1.$$ Hence, $\rep(Q,\beta)$ is a rational representation of the linearly reductive group $\GL(\beta)$ and the $\GL(\beta)-$orbits in $\rep(Q,\beta)$ are in one-to-one correspondence with the isomorphism classes of $\beta-$dimensional representations of $Q$. As $Q$ is a quiver without oriented cycles, one can show that there is only one closed $\GL(\beta)-$orbit in $\rep(Q,\beta)$ and hence the invariant ring $\text{I}(Q,\beta)= K [\rep(Q,\beta)]^{\GL(\beta)}$ is exactly the base field $K$.

Now, consider the subgroup $\SL(\beta) \subseteq \GL(\beta)$ defined by
$$
\SL(\beta)=\prod_{i \in Q_0}\SL_{\beta(i)}(K).
$$

Although there are only constant $\GL(\beta)-$invariant polynomial functions on $\rep(Q,\beta)$, the action of $\SL(\beta)$ on
$\rep(Q,\beta)$ provides us with a highly non-trivial ring of semi-invariants. Note that any $\sigma \in \ZZ^{Q_0}$ defines a
rational character of $\GL(\beta)$ by $$(g(i))_{i \in Q_0} \in \GL(\beta) \mapsto \prod_{i \in Q_0}(\det g(i))^{\sigma(i)}.$$ In this way, we can identify $\Gamma:=\ZZ ^{Q_0}$ with the group $X^\star(\GL(\beta))$ of rational characters of $\GL(\beta)$, assuming that $\beta$ is a sincere dimension vector. In general, we have only the natural epimorphism $\Gamma \to X^*(\GL(\beta))$. We also refer to the rational characters of $\GL(\beta)$ as (integral) weights.

Let $\SI(Q,\beta)= K[\rep(Q,\beta)]^{\SL(\beta)}$ be the ring of semi-invariants. As $\SL(\beta)$ is the commutator subgroup of
$\GL(\beta)$ and $\GL(\beta)$ is linearly reductive, we have
$$\SI(Q,\beta)=\bigoplus_{\sigma \in X^\star(\GL(\beta))}\SI(Q,\beta)_{\sigma},
$$
where $$\SI(Q,\beta)_{\sigma}=\lbrace f \in K[\rep(Q,\beta)] \mid g\cdot f= \sigma(g)f \text{~for all~}g \in \GL(\beta)\rbrace$$ is the space of semi-invariants of weight $\sigma$.

In a seminal paper \cite{K}, King constructed, via GIT, moduli spaces for finite-dimensional algebras. In what follows, we recall
King's main results. Note that the one-dimensional torus
$$
T=\{(t\Id_{\beta(i)})_{i \in Q_0} \mid t \in K^{*} \} \subseteq
\GL(\beta)
$$
acts trivially on $\rep(Q,\beta)$ and so there is a well-defined action of $\PGL(\beta)={\GL(\beta)/T}$ on $\rep(Q,\beta)$.

\begin{definition} \cite[Definition 2.1]{K} Let $\beta$ be a dimension vector of $Q$ and
$\sigma \in \ZZ^{Q_0}$ an integral weight. A representation $W \in \rep(Q,\beta)$ is said
to be:
\begin{enumerate}
\renewcommand{\theenumi}{\arabic{enumi}}

\item \emph{$\sigma$-semi-stable} if there exists a semi-invariant $f \in \SI(Q,\beta)_{m\sigma}$ with $m \geq 1$, such that $f(W)\neq 0$;

\item \emph{$\sigma$-stable} if there exists a semi-invariant $f \in \SI(Q,\beta)_{m\sigma}$ with $m \geq 1$, such that $f(W)\neq 0$ and, furthermore,
the $\GL(\beta)$-action on the principal open subset defined by $f$ is closed and $\dim \GL(\beta)W=\dim \PGL(\beta)$.
\end{enumerate}
\end{definition}
Note that any $\sigma$-stable representation is, in particular, a Schur representation. Consider the (possibly empty) open subsets
$$\rep(Q,\beta)^{ss}_{\sigma}=\{W \in \rep(Q,\beta)\mid W \text{~is~}
\text{$\sigma$-semi-stable}\}$$
and $$\rep(Q,\beta)^s_{\sigma}=\{W \in \rep(Q,\beta)\mid W \text{~is~}
\text{$\sigma$-stable}\}$$
of $\beta$-dimensional $\sigma$(-semi)-stable representations.

The GIT-quotient of $\rep(Q,\beta)^{ss}_{\sigma}$ by $\PGL(\beta)$ is
$$
\mathcal{M}(Q,\beta)^{ss}_{\sigma}=\Proj(\bigoplus_{m \geq 0}\SI(Q,\beta)_{m\sigma}).
$$
This is an irreducible projective variety whose closed points parameterize the closed $\GL(\beta)$-orbits
in $\rep(Q,\beta)^{ss}_{\sigma}$. For given $\beta, \sigma \in \RR^{Q_0}$, we
define
$$
\sigma(\beta)=\sum_{i \in Q_0}\sigma(i)\beta(i).
$$

In \cite{K}, King found a representation-theoretic description of the (semi-)stable representations and of the closed
orbits in $\rep(Q,\beta)^{ss}_{\sigma}$:

\begin{proposition}\cite[Proposition 3.1, 3.2]{K}\label{King-criteria} Let $\beta$ be a (non-zero) dimension vector and $\sigma$ an integral weight of $Q$. For a given representation
$W \in \rep(Q,\beta)$, the following are true:

\begin{enumerate}
\renewcommand{\theenumi}{\arabic{enumi}}

\item $W$ is $\sigma$-semi-stable if and only if $\sigma(\dd{W})=0$ and $\sigma(\dd{W'})
\leq 0$ for every subrepresentation $W'$ of $W$;

\item $W$ is $\sigma$-stable if and only if $\sigma(\dd{W})=0$ and $\sigma(\dd{W'})<0$ for
every proper subrepresentation $0 \neq W' \varsubsetneq W$;

\item $\GL(\beta)W$ is closed in $\rep(Q,\beta)^{ss}_{\sigma}$ if and only if $W$ is a direct sum of $\sigma$-stable representations. We call such a representation \emph{$\sigma$-poly-stable}.
\end{enumerate}
\end{proposition}
Note that we can use this result to define $\sigma$-(semi-)stable representations with respect to any real-valued function $\sigma \in \RR^{Q_0}$. We say that a dimension vector $\beta$ is \emph{$\sigma$(-semi)-stable} if there exists $\sigma$(-semi)-stable representation $W \in \rep(Q,\beta)$.

\subsection{The $\sigma$-stable decomposition} In this section, we recall Derksen and Weyman's \cite{DW2} notion of $\sigma$-stable decomposition of dimension vectors which proves to be a very powerful tool for studying semi-invariants of quivers.

Given a rational-valued function $\sigma \in \QQ^{Q_0}$, we define $\rep(Q)^{ss}_{\sigma}$ to be the full subcategory of $\rep(Q)$
consisting of all $\sigma$-semi-stable representations, i.e., those representations that satisfy Proposition \ref{King-criteria}{(1)}. Similarly, we define $\rep(Q)^{s}_{\sigma}$ to be the full subcategory of $\rep(Q)$ consisting of all $\sigma$-stable representations. (Of course, the zero representation is always semi-stable but not stable.)

It is easy to see that $\rep(Q)^{ss}_{\sigma}$ is a full exact subcategory, closed under extensions, and whose simple objects are precisely the $\sigma$-stable representations. Moreover, $\rep(Q)^{ss}_{\sigma}$ is Artinian and Noetherian, and hence, every $\sigma$-semi-stable representation has a Jordan-H{\"o}lder filtration in $\rep(Q)^{ss}_{\sigma}$.

Let $\alpha, \beta$ be two dimension vectors. We define
$$
\ext_Q(\alpha,\beta)=\min \{\dim_K \Ext^1_{Q}(V,W) \mid (V,W) \in \rep(Q,\alpha)\times \rep(Q,\beta)\}
$$
and
$$
\hom_Q(\alpha,\beta)=\min \{\dim_K \Hom_{Q}(V,W) \mid (V,W) \in \rep(Q,\alpha)\times \rep(Q,\beta)\}.
$$
It is not difficult to show that the dimensions of $\Ext^1_Q$ and $\Hom_Q$ spaces are upper-semi-continuous as functions on $\rep(Q,\alpha)\times \rep(Q,\beta)$. Hence, the above minimal values are attained on open subsets of $\rep(Q,\alpha) \times \rep(Q,\beta)$.

Let $\beta$ be a (non-zero) $\sigma$-semi-stable dimension vector where $\sigma \in \ZZ^{Q_0}$. We say that $$\beta=\beta_1 \pp \beta_2 \pp \ldots \pp \beta_s$$ is the \emph{$\sigma$-stable decomposition} of $\beta$ if a general representation in $\rep(Q,\beta)$ has a Jordan-H{\"o}lder filtration in $\rep(Q)^{ss}_{\sigma}$ with factors of dimensions $\beta_1, \ldots ,\beta_s$ (in some order). We write $c\cdot\beta$ instead of $\beta \pp \beta \pp \ldots \pp \beta$ ($c$ times).

The next proposition gives some basic properties of the dimension vectors occurring in the $\sigma$-stable decomposition of a
dimension vector. It is essential for proving Proposition \ref{prop-effective-wts}, Schofield's Embedding Theorem \ref{S-embed-thm}, and Theorem \ref{domains-exceptional-thm}.

\begin{proposition}\cite[Proposition 3.18]{DW2}\label{sigma-prop}
Let $\beta$ be a $\sigma$-semi-stable dimension vector and let
$$\beta=c_1\cdot \beta_1 \pp c_2\cdot \beta_2
\pp \ldots \pp c_l \cdot \beta_l$$
be the $\sigma$-stable decomposition of $\beta$ with the $\beta_i$ distinct. Then:
\begin{enumerate}
\renewcommand{\theenumi}{\arabic{enumi}}
\item the $\beta_i$ are Schur roots;

\item $\hom_Q(\beta_i,\beta_j)=0$ for all $i \neq j$;

\item after rearranging, we can assume that $\ext_Q(\beta_i,\beta_j)=0$ for all $1 \leq i<j \leq l$.
\end{enumerate}
\end{proposition}

\subsection{Domains of semi-invariants}
Let $\alpha$ and $\beta$ be two dimension vectors. We write $\alpha \hookrightarrow \beta$ if every representation of dimension vector $\beta$ has a subrepresentation of dimension vector $\alpha$.

Recall that if $\beta$ is a dimension vector of $Q$, the domain of semi-invariants associated to $(Q,\beta)$ is
$$
D(\beta)=\{ \alpha \in \QQ^{Q_0}  \mid \langle \alpha, \beta \rangle
=0 \text{~and~} \langle \alpha, \beta' \rangle \leq 0 \text{~for
all~} \beta' \hookrightarrow \beta \}.
$$

\begin{remark}\label{IOTW-rmk-sincere} Let $\beta$ be a dimension vector and $i \in Q_0$. Then, it is easy to see that $\beta_{i}=0$ if and only if $\gamma_i \in D(\beta)$ if and only if $-\gamma_i \in D(\beta)$ (see for example \cite[Lemma 6.5.6]{IOTW}).
\end{remark}

\begin{lemma}\cite[Lemma 6.5.7]{IOTW} \label{sincere-lemma} Let $\alpha, \beta \in \ZZ^{Q_0}$ be two integer valued functions.
\begin{enumerate}
\renewcommand{\theenumi}{\arabic{enumi}}

\item Assume that $\beta$ is a sincere dimension vector and $\alpha \in D(\beta)$. Then $\alpha$ is also a dimension vector.

\item Dually, if $\alpha$ is a sincere dimension vector, $\langle \alpha, \beta \rangle =0$, and $\langle \alpha',\beta \rangle \geq 0$ for
all $\alpha' \hookrightarrow \alpha$ then $\beta$ is also a dimension vector.
\end{enumerate}
\end{lemma}

\begin{remark} Note that the lemma above can also be deduced from \cite[Theorem 1]{DW1}.
\end{remark}

When $\beta$ is a sincere dimension vector, a description of the lattice points of $D(\beta)$ in terms of perpendicular categories was obtained independently in \cite{CB}, \cite{DW1}, and \cite{SVB}. An extension of this result to the case of arbitrary dimension vectors was obtained by Igusa-Orr-Todorov-Weyman in \cite{IOTW}. For a dimension vector $\delta$, we define
$$
P_{\delta}=\bigoplus_{i\in Q_0}P_i^{\delta(i)}.
$$

Now, we can state:

\begin{theorem} \cite{IOTW} \label{gen-Sat-thm} Let $\beta$ be a dimension vector of $Q$ and $\alpha \in \ZZ^{Q_0}$ an
integral weight.

\begin{enumerate}
\renewcommand{\theenumi}{\arabic{enumi}}
\item There are unique dimension vectors $\alpha'$ and $\delta$ such that
$$\alpha=\alpha'-\dd{P_{\delta}} \text{~and~} \supp(\alpha')\cap \supp(\delta)=\emptyset.$$
Furthermore, in the special case when $\alpha \in \ZZ^{Q_0}_{\geq 0}$, one has $\alpha=\alpha'$ and $\delta=0$.

\item The following statements are equivalent:
\begin{enumerate}
\renewcommand{\theenumi}{\arabic{enumi}}
\item $\alpha \in D(\beta)$;

\item there is an $\alpha'$-dimensional representation $V$ such that

\begin{enumerate}
\renewcommand{\theenumi}{\roman{enumi}}
\item $\Hom_{Q}(V,W)=\Ext^1_Q(V,W)=0$ for some (equivalently, a generic) $W \in \rep(Q,\beta)$;
\item $\supp(\beta)\cap \supp(\delta)=\emptyset$.
\end{enumerate}
\end{enumerate}
\end{enumerate}
\end{theorem}

\begin{proof} The first part of the theorem is proved in \cite[Lemma 5.3.2]{IOTW}. The second part follows from Proposition 5.1.4, Corollary 6.2.2 and Theorem 6.5.11 in \cite{IOTW}.
\end{proof}

If $\alpha \in \ZZ^{Q_0}$, we define the weight $\sigma=\langle \alpha,\cdot \rangle$ by
$$\sigma(i)=\langle \alpha,\varepsilon_i \rangle, \forall i \in
Q_0.$$ Conversely, it is easy to see that for any weight $\sigma \in \ZZ^{Q_0}$ there is a unique $\alpha \in \ZZ^{Q_0}$ (not
necessarily a dimension vector) such that $\sigma=\langle \alpha,\cdot \rangle$. Similarly, one can define $\mu = \langle \cdot,\alpha \rangle$.

\begin{proposition} \label{prop-effective-wts} Let $\beta$ be a dimension vector and $\sigma \in \ZZ^{Q_0}$ an integral weight.
\begin{enumerate}
\renewcommand{\theenumi}{\arabic{enumi}}
\item $\beta$ is $\sigma$-semi-stable if and only if $\sigma(\beta)=0$ and $\sigma(\beta')\leq 0$ for all $\beta' \hookrightarrow \beta$.

\item $\beta$ is $\sigma$-stable if and only if $\beta$ is non-zero, $\sigma(\beta)=0$, and $\sigma(\beta')<0$ for all $\beta' \hookrightarrow \beta$ and $\beta'\neq 0$ or $\beta$.
\end{enumerate}
\end{proposition}

\begin{remark} This result is undoubtedly well-known. The implication $"\Longrightarrow"$ of both $(1)$ and $(2)$ is proved in Proposition \ref{King-criteria}. The implication $"\Longleftarrow"$ of $(1)$ was proved independently in \cite{DW1} and \cite{SVB} (for a proof, see \cite[Theorem 2.4]{CB}) for the case where $\sigma$ is of the form $\sigma=\langle \alpha, \cdot \rangle$ with $\alpha$ a dimension vector. For the lack of a reference for arbitrary $\alpha$ or for the implication $"\Longleftarrow"$ of $(2)$, we include a proof below.
\end{remark}

\begin{proof} Working with the full subquiver of $Q$ whose set of vertices is $\supp(\beta)$ and using Lemma \ref{sincere-lemma}, we can assume that $\sigma=\langle \alpha, \cdot \rangle$ with $\alpha$ a dimension vector. The proof of $(1)$ now follows from the remark above.

Now, let us prove $"\Longleftarrow"$ of $(2)$. From $(1)$ we know that $\beta$ is $\sigma$-semi-stable and let us consider the $\sigma$-stable decomposition of $\beta$:
$$
\beta=c_1\cdot \beta_1 \pp c_2\cdot \beta_2
\pp \ldots \pp c_l \cdot \beta_l,
$$
where the $\beta_i$ satisfies the conditions $(1)-(3)$ of Proposition \ref{sigma-prop}. It is clear that $$\ext_Q(c_1\beta_1,\sum_{2\leq i \leq l}c_i\beta_i)=0$$ and hence
$c_1\beta_1 \hookrightarrow \beta$ by \cite[Theorem 3.2]{S1}. If $\beta$ is not $\sigma$-stable then either $\beta'=\beta_1$
(when $l=1$) or $\beta'=c_1\beta_1$ (when $l\geq 2$) is a proper dimension sub-vector of $\beta$ with $\beta' \hookrightarrow \beta$ and $\sigma(\beta')=0$. But this is a contradiction.
\end{proof}

\subsection{Derksen-Weyman Saturation and Reciprocity Properties}
Let $\alpha,\beta$ be two dimension vectors such that $\langle \alpha,\beta \rangle=0$. In \cite{S2}, Schofield discovered some very important semi-invariants of quivers. Consider the following polynomial function
$$
\begin{aligned}
c:\rep(Q,\alpha)&\times \rep(Q,\beta) \rightarrow K\\
c(V,W)&=\det(d^V_{W}).
\end{aligned}
$$
Note that $d^V_W$ from Ringel's canonical exact sequence (\ref{can-exact-seq}) is indeed a square matrix since $\langle \alpha,\beta \rangle=0$. Fix $(V,W)\in \rep(Q,\alpha) \times \rep(Q,\beta)$. Then it is easy to see that $c^V=c(V,\cdot):\rep(Q,\beta)\rightarrow K$ is a semi-invariant of weight $\langle \alpha, \cdot \rangle$ and $c_W=c(\cdot,W):\rep(Q,\alpha)\rightarrow K$ is a semi-invariant of weight $-\langle \cdot, \beta \rangle$.

\begin{remark}
We should point out that if $V$ is an $\alpha$-dimensional representation in $\rep(Q,\alpha)$, the semi-invariant $c^{V}$ is well-defined on $\rep(Q,\beta)$ up to a non-zero scalar.
\end{remark}

\begin{remark} Given $(V,W) \in \rep(Q,\alpha) \times \rep(Q,\beta)$, we have
$$
\Hom_Q(V,W)=\Ext_Q^1(V,W)=0 \Longleftrightarrow d^{V}_W \text{~is invertible~}
$$
and this implies that $V$ is $-\langle \cdot, \beta \rangle$-semi-stable and $W$ is $\langle \alpha, \cdot \rangle$-semi-stable.
\end{remark}

It is rather easy to see that the Schofield semi-invariants behave nicely with respect to exact sequences. In fact, we have:

\begin{lemma}\label{lemma-filt-inv}\cite[Lemma 1]{DW1} Let $\alpha$ and $\beta$ be two dimension
vectors such that $\langle \alpha,\beta \rangle =0$. Let $W$ be a $\beta$-dimensional representation which has a filtration
$$
F_{\bullet}(W): 0=W_0 \varsubsetneq W_1 \varsubsetneq \dots
\varsubsetneq W_{l-1} \varsubsetneq W_l=W,
$$
with $\langle \alpha, \dd{W_i / W_{i-1}}\rangle=0,
\forall 1 \leq i \leq l$. Then $ c_{W}=\prod_{i=1}^l c_{W_i /
W_{i-1}}$ on $\rep(Q,\alpha)$.
\end{lemma}

In \cite[Theorem 1]{DW1} (see also \cite{SVB} and \cite{DZ}), Derksen and Weyman proved a fundamental result showing that each weight space of semi-invariants is spanned by Schofield semi-invariants. This is known as the \emph{First Fundamental Theorem} for semi-invariants of quivers. Using the FFT, Derksen and Weyman derived some remarkable consequences.

\begin{theorem}[Saturation Theorem for Quivers] \cite[Theorem 2]{DW1} \label{Sat-thm} Let $\beta$ be a dimension vector and $\sigma \in \ZZ^{Q_0}$ a weight. Then the following are equivalent:
\begin{enumerate}
\renewcommand{\theenumi}{\roman{enumi}}
\item $\dim_K \SI(Q,\beta)_{\sigma}> 0$;
\item $\sigma(\beta)=0$ and $\sigma(\beta')\leq 0$ for all $\beta' \hookrightarrow \beta$.
\end{enumerate}
\end{theorem}

\begin{remark} It is worth pointing out that when $Q$ is the triple star quiver, the theorem above immediately implies the Saturation Conjecture for Littlewood-Richardson coefficients. See \cite{KT} and \cite[Corollary 2]{DW1}.
\end{remark}

We also have the so-called reciprocity property:

\begin{proposition}[Reciprocity Property] \cite[Corollary 1]{DW1} \label{reciproc-prop} Let $\alpha$ and $\beta$ be two dimension vectors. Then
$$\dim_K\SI(Q,\beta)_{\langle \alpha, \cdot \rangle} = \dim_K
\SI(Q,\alpha)_{-\langle \cdot, \beta \rangle}.$$
\end{proposition}

For two dimension vectors $\alpha$ and $\beta$, we define
$$
\alpha \circ \beta=\dim_K \SI(Q,\beta)_{\langle \alpha, \cdot \rangle} = \dim_K
\SI(Q,\alpha)_{-\langle \cdot, \beta \rangle}.
$$

The next lemma is especially useful for proving Theorem \ref{S-embed-thm} and Theorem \ref{domains-exceptional-thm}:

\begin{lemma}\label{rSchur-lemma} Assume that $\alpha$ is a dimension vector such that
$\GL(\alpha)$ acts with a dense orbit on $\rep(Q, \alpha)$ and let $\sigma$ be either $\langle \alpha, \cdot \rangle$ or $-\langle \cdot, \alpha \rangle$. If $\delta$ is a $\sigma$-stable dimension vector then $\delta$ is a real Schur root.
\end{lemma}

\begin{proof}
We know that the dimension vector of any stable representation is a Schur root. We only need to show that $\GL(\delta)$ acts with a dense orbit on $\rep(Q, \delta)$ which is equivalent to showing $$\dim_K \SI(Q,
\delta)_{\mu}\leq 1$$ for all weights $\mu$ of $\GL(\delta)$.

It is well-known (and easy to see) that $\GL(\alpha)$ acts with a dense orbit on $\rep(Q,\alpha)$ if and only if there exists $V \in \rep(Q,\alpha)$ with $\Ext^1_Q(V,V)=0$. Consequently, $\GL(n\alpha)$ acts with a dense orbit on $\rep(Q,n\alpha)$ for any integer $n>0$. From this and the Reciprocity Property \ref{reciproc-prop}, we deduce that
$$
\dim_K \SI(Q, \delta)_{n \sigma}=1,
$$
for any integer $n>0$.

Now, let $\mu$ be a weight such that $\SI(Q,\delta)_{\mu}\neq 0$; in particular, $\delta$ is $\mu$-semi-stable. As $\delta$ is $\sigma$-stable and using Proposition \ref{prop-effective-wts}, we can always find a sufficiently large integer $n>0$ such that $n\sigma(\delta')-\mu(\delta') \leq 0$ for all $\delta' \hookrightarrow \delta$. But this is equivalent to $\SI(Q, \delta)_{n \sigma-\mu} \neq \{0\}$ by the Saturation Theorem \ref{Sat-thm}. Multiplying the semi-invariants in $\SI(Q, \delta)_{\mu}$ by a fixed non-zero semi-invariant in $\SI(Q, \delta)_{n\sigma-\mu}$, we get an injective linear map from $\SI(Q, \delta)_{\mu}$ into $\SI(Q,\delta)_{n\sigma}$. Consequently, $\dim_K \SI(Q, \delta)_{\mu} \leq \dim_K \SI(Q, \delta)_{n\sigma}=1$, and so, $\delta$ is a real Schur root.
\end{proof}

\section{Schofield's Embedding Theorem} \label{Embed.sec}

In this section, we review Schofield's results on perpendicular categories from \cite{S2}. We give new proofs of his results by using some of the tools we have already discussed and the $A_{\infty}$-formalism.

For a given representation $V$, the right perpendicular category of $V$, denoted by $V ^\perp {}$, is the full subcategory of representations $W$ such that $\Hom_Q(V,W)=\Ext^1_Q(V,W)=0$ (we also write $V \perp W$ in this case). Similarly, one defines the left perpendicular category ${}^\perp V$.

Now, let $\alpha$ be a dimension vector. We define $\alpha ^\perp {}$ to be the full subcategory consisting of all those representations $W$ with $V \perp W$ for some (or equivalently generic) $V \in \rep(Q,\alpha)$. Similarly, we define ${}^\perp \alpha$.

A dimension vector $\alpha$ is said to be \emph{pre-homogeneous} if $\GL(\alpha)$ acts with a dense orbit on the representation space $\rep(Q,\alpha)$. Our goal in this section is to understand the categories $\alpha ^\perp {}$ (and ${}^\perp \alpha$) when $\alpha$ is a pre-homogeneous dimension vector. For such a dimension vector $\alpha$, we claim that
$$\alpha ^\perp {}=\rep(Q)^{ss}_\sigma \text{~and~} {}^\perp \alpha=\rep(Q)^{ss}_{\mu}$$
where $\sigma= \langle \alpha, \cdot \rangle$ and $\mu=-\langle \cdot, \alpha \rangle$. Indeed, if $V \in \rep(Q, \alpha)$ is a rigid representation then each non-zero weight space of semi-invariants of the form $\SI(Q,\beta)_{m\sigma}$ is one-dimensional being spanned by $(c^V)^m$. Using this observation it is now easy to see that $\alpha ^\perp {}=\rep(Q)^{ss}_\sigma$. Similarly, one can prove the other claim.

A sequence $\mathcal{E}=(\beta_1, \dots, \beta_l)$ of dimension vectors is said to be a \emph{quiver exceptional sequence} if
\begin{enumerate}
\renewcommand{\theenumi}{\roman{enumi}}
\item each $\beta_i$ is a real Schur root;
\item $\beta_i \perp \beta_j$ for all $1 \leq i<j \leq l$;
\item $\langle \beta_j, \beta_i \rangle \leq 0$ for all $1 \leq i<j \leq l$.
\end{enumerate}
If we drop condition (iii), we call $\mathcal{E}$ just an \emph{exceptional sequence}. A sequence $(E_1, \dots, E_l)$ of exceptional representations is said to be a (quiver) exceptional sequence (of representations) if $(\dd{E_1}, \dots,\dd{E_l})$ is a (quiver) exceptional sequence. We say that $\mathcal{E}$ is \emph{complete} if $l=|Q_0|$.

\begin{remark} \label{indep-stabledecom-rmk} Let $\sigma$ be an integral weight and $\beta$ a non-zero $\sigma$-semi-stable dimension vector. It follows from  Proposition \ref{sigma-prop} that the Schur roots occurring in the $\sigma$-stable decomposition of $\beta$ form, possible after reordering, a quiver exceptional sequence.
\end{remark}

\begin{remark} \label{quiverexceptional-rmk} Let $(\beta_1, \dots, \beta_l)$ be a quiver exceptional sequence. We claim that $\hom_Q(\beta_i,\beta_j)$ is zero for all $i \neq j$. This is clearly true for $i<j$. Since $\beta_i$ and $\beta_j$ are Schur roots and $\ext_Q(\beta_i,\beta_j)=0$ for $i<j$, we know that either $\hom_{Q}(\beta_j,\beta_i)=0$ or $\ext_{Q}(\beta_j,\beta_i)=0$ by \cite[Theorem 4.1]{S1}. From this and the fact that $\langle \beta_j,\beta_i \rangle \leq 0$, we finally deduce that $\hom_Q(\beta_j,\beta_i)=0$ for $i<j$. In particular, $\ext_Q(\beta_i,\beta_j)=-\langle \beta_i,\beta_j \rangle$ for all $1 \leq i \neq j \leq l$. Moreover, the matrix $(\langle \beta_i, \beta_j \rangle)_{i,j}$ is lower triangular with ones on the diagonal, and hence, the $\beta_i$ are linearly independent over $\RR$.
\end{remark}

Now, we are ready to give a new proof of Schofield's Embedding Theorem from \cite[Theorem 2.5]{S2}.

\begin{proof}[Proof of Theorem \ref{S-embed-thm}] We prove the theorem for $\sigma=\langle \alpha, \cdot \rangle$. The case where $\sigma=-\langle \cdot, \alpha \rangle$ is completely analogous.

(1) We know from Lemma \ref{rSchur-lemma} that every $\sigma$-stable representation is exceptional and that in each dimension vector there is at most one $\sigma$-stable representation up to isomorphism. Next, we claim that if $E_1, \dots, E_m$ are pairwise non-isomorphic $\sigma$-stable representations, their dimension vectors $\dd{E_i}$ must be linearly independent over $\ZZ$. Assume to the contrary that there are integers $k_l \in \ZZ_{>0}$ such that
$$
\beta=\sum_{i \in I} k_i \dd{E_i}=\sum_{j \in J} k_j \dd{E_j},
$$
with $I\bigcap J=\emptyset$. Set $E_I=\oplus_{i\in I}E_i$ and $E_J=\oplus_{j \in J}E_j$ and note that $E_I$ and $E_J$ are $\sigma$-poly-stable representations of the same dimension $\beta$. Since $\alpha$ is pre-homogeneous, we know that any of its positive integer multiples is pre-homogeneous, and hence, $\dim_K \SI(Q,l \alpha)_{-\langle \cdot, \beta \rangle}=1$ for all integers $l \geq 0$. By the Reciprocity Property \ref{reciproc-prop}, this is equivalent to $\dim_K \SI(Q,\beta)_{l\sigma}=1$, and so, the moduli space $\mathcal{M}(Q,\beta)^{ss}_{\sigma}$ is just a point. From Proposition \ref{King-criteria}{(3)} it follows that $E_I\cong E_J$ which is a contradiction. The first part of the theorem now follows.

(2) Since the $\beta_i$ are linearly independent, we deduce that $$\beta_0=\sum_{i=1}^l \beta_i$$ is the $\sigma$-stable decomposition of $\beta_0$. It follows from Proposition \ref{sigma-prop} that after rearranging $\mathcal{E}=(\beta_1, \dots, \beta_l)$ is a quiver exceptional sequence.

(3) Let $\filt(\mathcal{E})$ be the full subcategory of $\rep(Q)$ whose objects have a finite filtration with factors among the $E_i$. We clearly have that $\filt(\mathcal{E})=\rep(Q)^{ss}_{\sigma}$. Using the $A_{\infty}$-formalism, Keller \cite[Section 2.3]{Kel2} (see also \cite[Section 7.7]{Kel1}) proved that $\filt(\mathcal{E})$ is determined by the Yoneda algebra $\Ext^*_Q(\bigoplus_{i=1}^lE_i,\bigoplus_{i=1}^lE_i)$ equipped with its $A_{\infty}$-algebra structure. More precisely, let $\mathcal A$ be the $A_{\infty}$-category with objects $X_1, \dots, X_l$ and morphism spaces $\Hom^*_{\mathcal A}(X_i,X_j)=\Ext^*_Q(E_i,E_j)$ and let $\twist(\mathcal A)$ be the category of twisted stalks over $\mathcal A$. Since the are no higher $\Ext_Q^i$ spaces (with $i \geq 2$) over path algebras, the objects of $\twist(\mathcal A)$ can be described as pairs $(X, \delta)$, where $X=(\{X_i\}_{1 \leq i \leq l}, \{V_i\}_{1 \leq i \leq l})$, formally written as $X=\bigoplus_{i=1}^l V_i \otimes X_i$ with the $V_i$ finite dimensional vector spaces, called multiplicity spaces, and $\delta=(\delta_{ji})_{1 \leq j,i \leq l}$ is a matrix of morphisms $\delta_{ji} \in \Hom_K(V_i,V_j)\otimes \Hom^1_{\mathcal A}(X_i,X_j)$. (Note that $\delta$ is an upper triangular matrix and it satisfies the Maurer-Cartan equation.) Using the fact that $\Hom^0_{\mathcal A}(X_i,X_j)=\{0\}$ for $i \neq j$ and that it is just the base field $K$ when $i=j$, it is easy to see that $\twist(\mathcal A)=\rep(Q(\mathcal{E}))$. From \cite[Proposition 2.3]{Kel1} we get the desired equivalence of categories. The fact that the map $I$ preserves the Euler forms of $Q$ and $Q(\mathcal{E})$ follows immediately from formula $(\ref{Euler-formula})$.
\end{proof}

\begin{remark}\label{S-embed-rmk1} It follows from the proof above that $l$ is at most $|Q_0|$ minus the number $r$ of non-isomorphic indecomposable direct summands of the generic $\alpha$-dimensional representation. In fact, Schofield showed in \cite{S2} that $l=|Q_0|-r$.
\end{remark}

\begin{remark} \label{S-embed-rmk2} Let $\mathcal{E}=(E_1, \dots, E_l)$ be a quiver exceptional sequence and let $\filt(\mathcal{E})$ be the full subcategory of $\rep(Q)$ whose objects have a filtration with factors among the $E_i$. It is now clear that Theorem \ref{S-embed-thm}{(3)} remains true for $\mathcal{E}$ (see also \cite[Theorem 2.39]{DW2}).
\end{remark}

The next theorem, which is the main result of \cite[Theorem 4.3]{S2}, provides us with algebraically independent generators of the algebra of semi-invariants $\SI(Q,\alpha)$ for the case where $\alpha$ is pre-homogeneous. Although it is not needed for our direct purposes, we include a new proof for completeness:

\begin{theorem} \label{S-preh-semi-inv-thm} Let $\alpha$ be a sincere pre-homogeneous dimension vector. If $E_1, \dots, E_l$ are the pairwise non-isomorphic $\langle \alpha, \cdot \rangle$-stable representations then the Schofield semi-invariants $c_{E_i}$ are algebraically independent and
$$
\SI(Q,\alpha)=K[c_{E_1}, \dots, c_{E_l}].
$$
\end{theorem}

\begin{remark} Schofield's original proof of the theorem above uses Luna's \cite{L} \'{e}tale slice machinery and Sato-Kimura \cite{SK} classification of pre-homogeneous vector spaces. Our proof is much simpler in nature and it follows immediately from Theorem \ref{S-embed-thm} and some of the basic properties of the Schofield semi-invariants.
\end{remark}

\begin{proof} Note that the weights $-\langle \cdot, \dd{E_i} \rangle$ of the semi-invariants $c_{E_i}$ are linearly independent over $\ZZ$ by Theorem \ref{S-embed-thm}{(2)} and Remark \ref{quiverexceptional-rmk}. Therefore, these semi-invariants are algebraically independent. (To conclude this, we need the assumption that $\alpha$ is sincere.)

It remains to show that each non-zero weight-space $\SI(Q,\alpha)_{\mu}$ is spanned by a monomial in the $c_{E_i}$. Since $\alpha$ is sincere, we know that $\mu=-\langle \cdot, \beta \rangle$ with $\beta$ a dimension vector by Lemma \ref{sincere-lemma}. Now, it easy to see that $\SI(Q,\alpha)_{\mu}$ is spanned by a semi-invariant of the form $c_{W}$ with $W$ a $\sigma$-semi-stable $\beta$-dimensional representation where $\sigma=\langle \alpha, \cdot \rangle$. Consider a Jordan-H\"{o}lder filtration of $W$
$$
F_{\bullet}(W): 0=W_0 \varsubsetneq W_1 \varsubsetneq \dots
\varsubsetneq W_{l-1} \varsubsetneq W_l=W,
$$
with $W_i/W_{i-1}$ one of the $\sigma$-stable representations $E_j$. Using Lemma \ref{lemma-filt-inv}, we can write
$$
c_W=\prod_{i=1}^l c_{W_i/
W_{i-1}},
$$
and this finishes the proof.
\end{proof}

\section{Cluster fans and cones of finite-stability conditions} \label{cluster-fan-stability.sec} In this section, we first recall the construction of the cluster fan of a quiver $Q$ without oriented cycles. Recall that the set of almost positive real Schur roots is
$$\Psi(Q)_{\geq -1}=\{\beta \mid \beta \text{~is a real Schur root of $Q$~}\} \cup \{-\gamma_i \mid i \in Q_0\}.$$
We should point out that in general the set of all real Schur roots depends on the orientation of $Q$.

To construct the (possibly infinite) cluster fan $C(Q)$ on the ground set $\Psi(Q)_{\geq -1}$, we need some definitions first. If $\beta_1, \beta_2 \in \Psi(Q)_{\geq -1}$, their compatibility degree is defined by
$$
(\beta_1 || \beta_2)_Q=
\begin{cases}
\ext_Q(\beta_1,\beta_2)+\ext_Q(\beta_2,\beta_1)& \text{if the
$\beta_i$ are real Schur roots}\\
\beta(i) & \text{if $\{\beta_1,\beta_2\}=\{\beta,-\gamma_i\}$
with $\beta$ a real Schur root}\\
0 & \text{otherwise}
\end{cases}
$$

A subset $C \subseteq \Psi(Q)_{\geq -1}$ is said to be compatible if $(\beta_1 || \beta_2)_Q=0$ for all $\beta_1, \beta_2 \in C$. A maximal (with respect to inclusion) compatible set is called a \emph{cluster}.

\begin{remark} Note that the compatibility degrees are dimensions of $\Ext$ spaces between indecomposable objects in the cluster category $\mathcal C_Q$ associated to $Q$ (see \cite{BMRRT}). Now, let $\phi$ be the bijection that sends $-\gamma_i$ to $-\varepsilon_i$ and is the identity map on the set of real Schur roots. Then $C \subseteq \Psi(Q)_{\geq -1}$ is compatible if and only if $\phi(C)$ is compatible in the sense of \cite{MMZ}, assuming $Q$ is a Dynkin quiver.
\end{remark}

\begin{remark} It is not difficult to see that any compatible subset $C$ of $\Psi(Q)_{\geq -1}$ is linearly independent over $\RR$. Indeed, write $C=\{-\gamma_{i_1}, \dots, \gamma_{i_l}, \beta_{l+1}, \dots, \beta_n \}$ and assume, without loss of generality, that $\gamma_{i_j} \perp \gamma_{i_k}$ for all $1 \leq j<k \leq l$. (This is always possible since $Q$ has no oriented cycles.) Next, using \cite[Theorem 2.4]{S1}, we can rearrange the $\beta_j$ so that $\hom_Q(\beta_m,\beta_p)=0$ for all $l+1 \leq m<p \leq n$. Now, let $\alpha_j=\gamma_{i_j}$ for $1 \leq j \leq l$, and $\alpha_j=\beta_j$ for $l+1 \leq j \leq n$. Note that the matrix $(\langle \alpha_i,\alpha_j \rangle)_{i,j}$ is lower triangular with ones on the diagonal, and hence, is invertible. Consequently, the elements of $C$ are linearly independent over $\RR$; in particular, $C$ must have at most $|Q_0|$ linearly independent elements.
\end{remark}

For $C\subseteq \QQ^N$ a finite set of points, let $\cone(C)$ be the rational convex polyhedral cone (in $\QQ^N$) generated by $C$, i.e., $\cone(C)=\{x \in \QQ^N \mid x=\sum_{c \in C}\lambda_c c \text{~with~} \lambda_c \in \QQ_{\geq 0} \}$. Next, let us record the following well-known result which can be easily proved using basic results from tilting theory:

\begin{theorem}\label{cluster-fan-thm} Let $Q$ be a quiver with $N$ vertices. Then the collection of cones
$$
C(Q)=\{\cone(C) \mid C \text{~is a compatible subset of~} \Psi(Q)_{\geq -1} \}
$$
is a smooth fan of pure dimension $N$. We call $C(Q)$ the \emph{cluster fan} of $Q$.
\end{theorem}

Our goal is to give a more geometric interpretation of the cluster fan of $Q$. Recall that the cone (not necessarily convex) of finite stability conditions is the set of all $\sigma \in \QQ^{Q_0}$ for which there are finitely many (possibly none) $\sigma$-stable representations up to isomorphism.

\begin{proposition} \label{finite-stability-prop} Let $\sigma=\langle \alpha, \cdot \rangle \in \ZZ^{Q_0}$ be a weight with $\alpha \in \ZZ^{Q_0}.$
The following are equivalent:
\begin{enumerate}
\renewcommand{\theenumi}{\arabic{enumi}}
\item $\sigma \in \mathcal{S}(Q)$;

\item $\alpha \in \cone(C)$ for some compatible subset $C$ of $\Psi(Q)_{\geq -1}$.

\end{enumerate}
\end{proposition}

\begin{proof} First, let us prove the implication $"\Longrightarrow"$. Let $\beta_1, \dots, \beta_l$ be the $\sigma$-stable dimension vectors.
We clearly have that $\beta_0$ is $\sigma$-semi-stable where $\beta_0=\sum_{i=1}^l \beta_i$. (In case there are no $\sigma$-stable dimension vectors, we set $\beta_0=0$). From Theorem \ref{gen-Sat-thm}, we know
$$
\alpha=\alpha'-\dd{P_{\delta}},
$$
where $\supp(\alpha') \cap \supp(\delta)=\supp(\beta_0) \cap \supp(\delta)=\emptyset$, $\alpha' \in D(\beta_0)$. Let $Q'$ be the full sub-quiver of $Q$ with $Q'_0=\supp(\alpha')$ and denote by the same letter the restriction of $\alpha'$ to $Q'$.

Since there are only finitely many $\sigma$-stable representations, we know that there are only finitely many $\sigma$-polystable representations in each dimension vector. In other words, each moduli space $\mathcal M(Q,\beta)^{ss}_{\sigma}$ is either empty or a point and so
\begin{eqnarray}\label{eqn1}
\dim_K \SI(Q,\beta)_{\sigma} \leq 1,
\end{eqnarray}
for each dimension vector $\beta$.

Let $\beta'$ be a dimension vector of $Q'$ and extend it trivially to a dimension vector $\beta$ of $Q$. Denote the dimension vector of $P_{\delta}$ by $\alpha''$. From Remark \ref{IOTW-rmk-sincere} and the fact that $\beta$ and $\delta$ have disjoint supports, we deduce that $-\alpha'', \alpha'' \in D(\beta)$ which is equivalent to $\beta$ being semi-stable with respect to both $\langle \alpha'', \cdot \rangle$ and
$-\langle \alpha'', \cdot \rangle$ by Proposition \ref{prop-effective-wts}. From this observation and $(\ref{eqn1})$, one can easily see that $(\alpha' \circ \beta')_{Q'}=(\alpha' \circ \beta)_{Q} \leq 1$. As $\alpha'$ is a sincere dimension vector of $Q'$, an effective weight of $\GL(\alpha')$ is of the form $-\langle \cdot, \beta' \rangle$ with $\beta'$ some dimension vector of $Q'$ by Lemma \ref{sincere-lemma}. This shows that
$$
\dim_K \SI(Q', \alpha')_{\sigma'} \leq 1,
$$
for all weights $\sigma'$ of $Q'$, and consequently, $\GL(\alpha')$ acts with a dense orbit on $\rep(Q',\alpha')$. So, the dimension vectors of the indecomposable direct summands of the generic $\alpha'$-dimensional representation of $Q'$ form a compatible subset of $\Psi(Q')_{\geq -1}$. The (trivial) extension of this compatible subset to a compatible subset of $\Psi(Q)_{\geq -1}$ together with $\{-\gamma_i \mid i \in \supp(\delta)\}$ form a compatible subset $C \subseteq \Psi(Q)_{\geq -1}$ with $\alpha \in \cone{C}$.

To prove the other implication $"\Longleftarrow"$, let us assume that
$$
n\alpha=\sum_{j=1}^l \eta(j)\beta_j- \dd{P_{\delta}},
$$
where $n$ and the $\eta(j)>0$ are positive integers and $\{ \beta_1,\dots,\beta_l\}\bigcup\{-\gamma_i \mid i \in \supp(\delta)\}$ is a compatible subset of $\Psi(Q)_{\geq -1}$. Denote $\sum_{j=1}^{l}\eta(j)\beta_j$ by $\alpha'$ and $\dd P_{\delta}$ by $\alpha''$. Note that $\alpha'$ and $\delta$ have disjoint supports, and $\GL(\alpha')$ acts with a dense orbit on $\rep(Q,\alpha')$.

Now, let $\beta$ be a $\sigma$-stable dimension vector. Since this is equivalent to $\beta$ being $n\sigma$-stable, we can assume without loss of generality that $n=1$. From Theorem \ref{gen-Sat-thm}, we know that $\alpha' \in D(\beta)$ and $\supp(\beta) \cap \supp(\delta)=\emptyset$. In particular, we have $\langle \alpha', \beta' \rangle=\langle \alpha, \beta'
\rangle$ for all $\beta' \leq \beta$ (coordinatewise).

Using Proposition \ref{prop-effective-wts}, we deduce that $\beta$ is $\langle \alpha', \cdot \rangle$-stable. As $\GL(\alpha')$ acts with a dense orbit on $\rep(Q,\alpha')$ we know that there are only finitely many $\langle \alpha', \cdot \rangle$-stable dimension vectors by Theorem \ref{S-embed-thm}. This finishes the proof.
\end{proof}

\begin{remark}\label{rSchur-rmk} The last part of the proof above together with Lemma \ref{rSchur-lemma} shows that if $\beta$ is $\langle \alpha, \cdot \rangle$-stable with $\alpha \in \mathcal{S}(Q)$ then $\beta$ is a real Schur root.
\end{remark}

\begin{proof}[Proof of Theorem \ref{stabilitycone-thm}] It now follows from Theorem \ref{cluster-fan-thm} and Proposition \ref{finite-stability-prop}.
\end{proof}

The cone $\widetilde{\mathcal{S}}(Q)$ of \emph{effective finite stability conditions} of $Q$ is, by definition, the set of all $\sigma \in \QQ^{Q_0}$ for which there exists at least one, but finitely many $\sigma$-stable representations up to isomorphism.

\begin{theorem}\label{effectivestabilitycone-thm} Let $Q$ be a quiver with $N$ vertices. Then the cone $\widetilde{\mathcal{S}}(Q)$ has a fan covering given by the dual of the $(N-1)$-skeleton of the cluster fan of $Q$.
\end{theorem}

\begin{proof} It follows from Remark \ref{S-embed-rmk1} that $\alpha \in \rel(\cone(C))$ for some cluster $C$ if and only if there are no $\langle \alpha, \cdot \rangle$-stable representations. From this observation and Theorem \ref{stabilitycone-thm} we obtain the desired result.
\end{proof}

Now, we are ready to prove Theorem \ref{domains-cluster-thm}.

\begin{proof}[Proof of Theorem \ref{domains-cluster-thm}] From Theorem \ref{effectivestabilitycone-thm} and Remark \ref{rSchur-rmk}, we deduce that
\begin{equation}\label{incl-label}
\bigcup_{C} \cone(C)=\{\alpha \in \QQ^{Q_0} \mid \langle \alpha, \cdot \rangle \in \widetilde{\mathcal{S}}(Q) \} \subseteq \bigcup_{\beta} D(\beta),
\end{equation}
where the union on the left is over all compatible sets $C$ with (at most) $|Q_0|-1$ elements while the union of the right is over all real Schur roots $\beta$.

Now, let us prove the implication $``\Longrightarrow''$. First, let us look into the case where $Q$ is a Dynkin quiver. If $\alpha \in D(\beta)$ then there is at least one $\langle \alpha, \cdot \rangle$-stable representation, and furthermore, there can be only finitely many, up to isomorphism, stable representations as $Q$ is a Dynkin quiver. So, the inclusion in $(\ref{incl-label})$ is an equality for Dynkin quivers.

Next, let us assume that $Q$ is a generalized Kronecker quiver. Pick an  $\alpha \in D(\beta) \cap \ZZ^{Q_0}$ where $\beta$ is a real Schur root. From Theorem \ref{gen-Sat-thm}, we know that
$$
\alpha=\alpha'-\dd P_{\delta},
$$
where $\alpha'$ and $\delta$ are dimension vectors such that $\supp(\alpha') \cap \supp(\delta)=\supp(\beta) \cap \supp(\delta)=\emptyset$ and $\alpha' \in D(\beta)$.

If $\delta$ is the zero dimension vector then $\alpha=\alpha'$ is a $-\langle \cdot, \beta \rangle$-semi-stable dimension vector. Looking at the $-\langle \cdot, \beta \rangle$-stable decomposition of $\alpha$ and using Lemma \ref{rSchur-lemma}, we see that the Schur roots that occur in this decomposition of $\alpha$ are real Schur roots. Since the space of all vectors $\alpha'' \in \QQ^{Q_0}$ with $\langle \alpha'', \beta \rangle=0$ is one dimensional, we deduce that $\alpha$ is just a positive integer multiple of a real Schur root, i.e., $\alpha \in \cone(C)$ with $C$ a compatible subset with one element.

If $\delta$ is not the zero dimension vector then $\beta$ is just one of the two simple roots while $\alpha'$ must be the zero dimension vector. So, the inclusion in $(\ref{incl-label})$ is an equality for generalized Kronecker quivers, as well.

For the other implication $``\Longleftarrow''$, let $W$ be an exceptional representation and $V$ a representation such that $V \perp W$; in particular, $\dd V \in D(\dd W)$. It follows from Theorem \ref{gen-Sat-thm}{(1)} that
$$
n \dd V=\sum_{j=1}^l k_j\beta_j,
$$
where $n \geq 1$ is an integer, the $k_j$ are non-negative integers, and $\{\beta_1, \dots, \beta_l\}$ is a compatible subset of $\Psi(Q)_{\geq -1}$. This implies that $\langle \dd V, \dd V \rangle >0$ for all non-zero representations $V$ with $V \perp W$. Consequently, if $W$ is the indecomposable injective representation at some vertex $i$, we get that the quiver $Q\setminus \{i\}$ is a (union of) Dynkin quivers. Hence, $Q$ is either a Dynkin, or a generalized Kronecker quiver, or a Euclidean quiver with at least three vertices. In what follows, we show that the last case cannot occur.

Assume to the contrary that $Q$ is a Euclidean quiver with at least three vertices. Denote by $\delta_Q$ the isotropic Schur root of $Q$ and choose a vertex $i$ such that $Q \setminus \{i\}$ is a Dynkin quiver. Without loss of generality, let us assume that $i$ is a source. For $\beta_1=\delta_Q-\varepsilon_i$ and $\beta_2=\varepsilon_i$, we can see that $\mathcal{E}=(\beta_1,\beta_2)$ is a quiver exceptional sequence with $\langle \beta_2,\beta_1 \rangle=-2$. Hence, $Q(\mathcal{E})$ is the Kronecker quiver:
$$
\xy     (0,0)*{\cdot}="a";
        (10,0)*{\cdot}="b";
        {\ar@2{<-} "a";"b" };
\endxy
$$

Since $\mathcal{E}$ is not a complete exceptional sequence, we can always find a real Schur root $\beta$ such that $\beta_1, \beta_2 \in D(\beta)$. Indeed, this follows from the extension theorem for exceptional sequences due to Crawley-Boevey \cite{CB2}. Next, using Remark \ref{S-embed-rmk2}, we deduce that the Tits quadratic form of $Q(\mathcal{E})$ is weakly positive definite which is a contradiction.
\end{proof}

We end this section with some observations about the cluster fan and the GIT-classes of a quiver $Q$. Given two weights $\sigma_1, \sigma_2 \in \QQ^{Q_0}$, we say that they are \emph{GIT-equivalent} if $$\rep(Q)^{ss}_{\sigma_1}=\rep(Q)^{ss}_{\sigma_2}.$$ The \emph{GIT-class of a weight $\sigma \in \QQ^{Q_0}$}, denoted by $\langle \sigma \rangle$, is
$$
\langle \sigma \rangle=\{\sigma' \in \QQ^{Q_0} \mid \rep(Q)^{ss}_{\sigma}=\rep(Q)^{ss}_{\sigma'}\}.
$$

Now, let $C=\{\beta_j \mid 1 \leq j \leq l\} \cup \{-\gamma_{i_k} \mid l+1 \leq k \leq m \}$ be a compatible subset of $\Psi(Q)_{\geq -1}$ and pick $\alpha=\sum_{j=1}^l \eta(j)\beta_j-\sum_{k=l+1}^m c_k \gamma_{i_k} \in \cone(C)$ with $\eta(j)$ and $c_k$ positive integers. Denote $\langle \alpha, \cdot \rangle$ by $\sigma$. It follows from the proof of Theorem \ref{stabilitycone-thm} that
$$
\rep(Q)^{ss}_{\sigma}=\{W \in \rep(Q) \mid \beta_j \perp W, 1 \leq j \leq l, \text{~and~} W(i_k)=0, l+1 \leq k \leq m\}.
$$
(Here, by $\beta \perp W$, we simply mean that $U \perp W$ for some $\beta$-dimensional representation $U$.) Let $\alpha_C=\sum_{j=1}^l \beta_j-\sum_{k=l+1}^m \gamma_{i_k}$ and $\sigma_C=\langle \alpha_C, \cdot \rangle$. It is easy so see that
$$
\rel(I(\cone(C))) \subseteq \langle \sigma_C \rangle.
$$

It is clear that the inclusion above is strict whenever $C$ is a cluster. In fact, if $C$ is a cluster then $\rep(Q)^{ss}_{\sigma_C}$ consists of only the zero representation, and so, we have
$$
\langle \sigma_C \rangle=\bigcup_{C'} \rel(I(\cone(C'))),
$$
where the union on the right is over all clusters $C'$. That is to say, the clusters form one single GIT-class.

Now, let assume that $Q$ has at least three vertices and let $\beta$ be a real Schur root of $Q$. Using Theorem \ref{S-embed-thm}, we can always find compatible sets $C_1$ and $C_2$, each consisting of $|Q_0|-1$ real Schur roots of $Q$, such that the $\alpha_i:=\sum_{\alpha \in C_i} \alpha$, $1 \leq i \leq 2$, are two distinct pre-homogeneous dimension vectors for ${}^\perp \beta$. Denote $\langle \alpha_i, \cdot \rangle$ by $\sigma_i$, $1 \leq i \leq 2$. Using Theorem \ref{S-embed-thm} again, we see that $\sigma_1$ and $\sigma_2$ are GIT-equivalent since the $\sigma_i$-stable representations are precisely the $\beta$-dimensional exceptional representations for each $i \in \{1,2\}$. We conclude that
$$\rel(I(\cone(C_1))) \cup \rel(I(\cone(C_2))) \subseteq \langle \sigma_1 \rangle=\langle \sigma_2 \rangle.$$

So, the GIT-equivalence relation does not distinguish among the relative interiors of the cones generated by the compatible subsets of $\Psi(Q)_{\geq -1}$. Nonetheless, it would be interesting to find a (geometric) equivalence relation on $\mathcal{S}(Q)$ such that equivalence classes are precisely the relative interiors of the cones $I(\cone(C))$ with $C$ compatible subsets of $\Psi(Q)_{\geq -1}$.

\section{Domains of semi-invariants and quiver exceptional sets}\label{domains-exceptional.sec}
Our goal in this section is to find an extension of \cite[Theorem 8.1.7]{IOTW} to arbitrary quivers by keeping the domains of semi-invariants in our attention. For this, we need to work with (quiver) exceptional sets instead of compatible sets.

Let
$$\mathcal{E}=\{\beta_1, \dots,
\beta_l,-\gamma_{i_{l+1}}, \dots, -\gamma_{i_m} \}$$ be a subset of $\Psi(Q)_{\geq -1}$. We say that $\mathcal{E}$ is a
\emph{(quiver) exceptional set} if

\begin{enumerate}
\renewcommand{\theenumi}{\arabic{enumi}}
\item $\beta_j(i_k)=0$, $1 \leq j \leq l$, $l+1 \leq k \leq m$,

\item the $\beta_j$ can be rearranged so that $(\beta_1, \dots, \beta_l)$ is a (quiver) exceptional sequence.
\end{enumerate}

\begin{proof}[Proof of Theorem \ref{domains-exceptional-thm}] Let us denote $-\langle \cdot, \beta \rangle$ by $\sigma$. We know that there are only finitely many $\sigma$-stable dimension vectors and they are real Schur roots by Theorem \ref{S-embed-thm}. Let $\mathcal{E}_1, \dots, \mathcal{E}_m$ be all quiver exceptional sets such that each of them consists of $\sigma$-stable dimension vectors and integral vectors of the form $-\gamma_i$ with $i \notin \supp(\beta)$. By construction, we have
$$
\bigcup_{1 \leq i \leq m} \cone(\mathcal{E}_i)\subseteq D(\beta).
$$
Note that the size of each of the $\mathcal{E}_i$ is at most $|Q_0|-1$.

Now, let us prove the other inclusion. Pick $\alpha \in D(\beta) \bigcap \ZZ^{Q_0}$. From Theorem \ref{gen-Sat-thm}, we
know that there are dimension vectors $\alpha', \delta$ such that $\alpha'$ is $\sigma$-semi-stable, $\alpha=\alpha'-\dd P_{\delta}$, and $\supp(\beta) \cap \supp(\delta)=\supp(\alpha') \cap \supp(\delta)=\emptyset$. If $\supp(\delta)=\{i_{l+1} \dots, i_m\}$ then each $\gamma_{i_k}$ is in $D(\beta)$ as $\beta(i_k)=0$, and of course, $\dd P_{\delta}$ is a nonnegative linear combination of the $\gamma_{i_k}$.

Now, consider the $\sigma$-stable decomposition of $\alpha'$:
$$
\alpha'=c_1\cdot \beta_1 \pp \cdots \pp c_l \cdot \beta_l,$$
where the $\beta_i$ are distinct $\sigma$-stable dimension vectors and the $c_i$ are positive integers. From Lemma \ref{rSchur-lemma}, it follows that the $\beta_i$ are real Schur roots. Moreover, we know that after rearranging $(\beta_1, \dots, \beta_l)$ is a quiver exceptional sequence by Proposition \ref{sigma-prop}. Consequently, the set $\mathcal{E}=\{ \beta_1, \dots, \beta_l, -\gamma_{i_{l+1}}, \dots, -\gamma_{i_m}\}$ is one of the $\mathcal{E}_i$, and furthermore, $\alpha \in \cone(\mathcal{E}) \subseteq D(\beta)$. This finishes the first part of our theorem.

To prove the last part, let $\mathcal{E}=\{ \alpha_1, \dots, \alpha_l \}$ be a quiver exceptional set with $l \leq N-1$. If
$\alpha_k=-\gamma_{i_k}$ for all $1 \leq k \leq l$ then one can choose $\beta$ to be the simple root corresponding to some vertex
$i \in Q_0 \setminus \{i_1, \dots, i_l\}$. For such $\beta$, we clearly have $\mathcal{E} \subseteq D(\beta)$.

Now, let assume that
$$\mathcal{E}=\{\alpha_1=\beta_1, \dots, \alpha_l= \beta_l, \alpha_{l+1}=-\gamma_{i_{l+1}}, \dots, \alpha_m=-\gamma_{i_m}, \},$$ with $1 \leq l \leq m \leq n-1$. Then we can rearrange the $\beta_i$ so that $(\beta_1, \dots, \beta_l)$ is an exceptional sequence for
$\widetilde{Q}=Q \setminus \{i_{l+1}, \dots, i_m\}$. From the extension theorem for exceptional sequences due to Crawley-Boevey \cite{CB2}, we know that there exists a real Schur root $\beta$ of $\widetilde{Q}$ such that $\beta_i \in D_{\widetilde{Q}}(\beta)$. Extend $\beta$ (trivially) to a real Schur root of $Q$. Then, $\gamma_{i_k} \in D(\beta)$ as $i_k \notin \supp(\beta)$, and so, $\mathcal{E} \subseteq D(\beta)$.
\end{proof}

\begin{remark} We should point out that in case $\beta$ is a sincere dimension vector, it follows from Theorem \ref{S-embed-thm} and Lemma \ref{sincere-lemma} that $D(\beta)=\cone(\mathcal{E})$ where $\mathcal{E}$ is the quiver exceptional set consisting of all $-\langle \cdot, \beta \rangle$-stable dimension vectors. However, this fails when $\beta$ is not sincere. Indeed, if $\beta(i)=0$ for some $i \in Q_0$, the cone $D(\beta)$ is not strongly convex as it contains both $\gamma_i$ and $-\gamma_i$. So, $D(\beta)$ cannot even be simplicial in the non-sincere case.
\end{remark}

\begin{remark} We would like to point out that Theorem \ref{domains-exceptional-thm} remains true if instead of quiver exceptional sets we work with just exceptional sets.
\end{remark}

Let $\mathcal{E}=\{\alpha_1=\beta_1, \dots,\alpha_l= \beta_l, \alpha_{l+1}=-\gamma_{i_{l+1}}, \dots, \alpha_m=-\gamma_{i_m} \}$ be a quiver exceptional set. Define $Q(\mathcal E)$ to be the quiver with vertices $1, \dots, l$, and with $-\langle \beta_i,\beta_j \rangle$ arrows from vertex $i$ to vertex $j$ for all $1\leq i \neq j \leq l$. Recall that for a quiver exceptional set $\mathcal{E}$, $\ext(\beta_i,\beta_j)=-\langle \beta_i,\beta_j \rangle$ for all $1 \leq i \neq j \leq l$. For this reason, we also call $Q(\mathcal{E})$ the $\Ext$-quiver of $\mathcal{E}$.

We call a quiver exceptional set \emph{representation-finite} if $Q(\mathcal{E})$ is a (union of) Dynkin quivers.

\begin{example}\label{qes-ex} Let us give some examples of quiver exceptional sets whose $\Ext$-quivers are easy to describe.

\begin{itemize}
\item(Dynkin case) Let $Q$ be a Dynkin quiver and $\mathcal{E}$ a quiver exceptional set. From Theorem \ref{S-embed-thm}{(3)} and formula $(\ref{Euler-formula})$, we deduce that the Tits quadratic form of $Q(\mathcal{E})$ is weakly positive and hence $Q(\mathcal{E})$ is also a Dynkin quiver.

\item (Euclidean case) Let $Q$ be a Euclidean quiver and denote by $\delta_Q$ the isotropic Schur root of $Q$. Choose $i$ to be a vertex such that $Q \setminus \{i\}$ is a Dynkin quiver. Without loss of generality, let us assume that $i$ is a source. In this case, we take $\beta_1=\delta_Q-\varepsilon_i$ and $\beta_2=\varepsilon_i$. Then, the set $\mathcal{E}=\{\beta_1,\beta_2\}$ is a quiver exceptional set with $\langle \beta_2,\beta_1 \rangle=-2$. Hence, $Q(\mathcal{E})$ is the Kronecker quiver:
$$
\xy     (0,0)*{\cdot}="a";
        (10,0)*{\cdot}="b";
        {\ar@2{<-} "a";"b" };
\endxy
$$

\item (Wild case) Let $T_{4,3,4}$ be the wild star quiver with the following
orientation:
$$
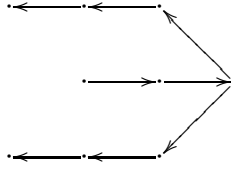
\xy (0, 0)*{\cdot}="a";
        (-10, 10)*{\cdot}="b";
        (-20,10)*{\cdot}="c";
        (-30,10)*{\cdot}="d";
        (-20,0)*{\cdot}="e";
        (-10,0)*{\cdot}="f";
        (-30,-10)*{\cdot}="g";
        (-20,-10)*{\cdot}="h";
        (-10,-10)*{\cdot}="i";
        {\ar@{->} "a";"b"};
        {\ar@{->} "b";"c"};
        {\ar@{->} "c";"d"};
        {\ar@{->} "f";"a"};
        {\ar@{->} "e";"f"};
        {\ar@{->} "a";"i"};
        {\ar@{->} "i";"h"};
        {\ar@{->} "h";"g"};
    \endxy
$$
Let us consider the exceptional set $\mathcal{E}=\{\beta_1, \beta_2\}$ of $T_{4,3,4}$
given by:
$$
\beta_1=
\begin{matrix}
1&2&3& \\
&0&3&4,\\
1&2&3&
\end{matrix}
$$
and
$$
\beta_2=
\begin{matrix}
0&0&0& \\
&1&0&0.\\
0&0&0&
\end{matrix}
$$

Since $\langle \beta_2, \beta_1 \rangle =-3$, $T_{4,3,4}(\mathcal{E})$ is the generalized Kronecker quiver:
$$
\xy     (0,0)*{\cdot}="a";
        (10,0)*{\cdot}="b";
        {\ar@3{<-} "a";"b" };
\endxy
$$
\end{itemize}
\end{example}

Next, we compare quiver exceptional sets with clusters:

\begin{proposition}\label{rep-finite-exceptional-prop} If $\mathcal E$ is a representation-finite quiver exceptional set then there are finitely many compatible sets $C_1, \dots, C_r$ such that
$$
\cone(\mathcal{E})=\bigcup_{i=1}^r \cone(C_i).
$$

\end{proposition}

\begin{proof} Write $\mathcal{E}=\{\alpha_1=\beta_1, \dots,\alpha_l= \beta_l, \alpha_{l+1}=-\gamma_{i_{l+1}}, \dots, \alpha_m=-\gamma_{i_m} \}$ and let
$$
\alpha=\sum_{j=1}^l\eta(j) \beta_j-\sum_{k=l+1}^m c_k\gamma_{i_k} \in \cone(\mathcal{E})
$$
with $\eta(j)$ and $c_k$ non-negative integers. We can assume that $(\beta_1, \dots,\beta_l)$ is a quiver exceptional sequence which, by some abuse, we denote by the same letter $\mathcal{E}$. From Theorem \ref{S-embed-thm}{(3)}, we know that there exists a full exact embedding of $\rep(Q(\mathcal{E}))$ into $\rep(Q)$ and let $I: \NN^{Q(\mathcal{E})_0}=\NN^l \rTo \NN^{Q_0}$ be the isometry induced by $\mathcal{E}$.

Since $Q(\mathcal{E})$ is a Dynkin quiver, we know that $\eta=(\eta(1), \dots, \eta(l))$ is a pre-homogeneous dimension vector. If $\eta_1, \dots, \eta_r$ are the dimension vectors of the indecomposable direct summands of a $\eta$-dimensional rigid representation of $Q(\mathcal{E})$ then it is easy to see that $$C=\{I(\eta_1), \dots, I(\eta_r)\} \bigcup \{-\gamma_{i_{l+1}}, \dots, -\gamma_{i_m} \}$$ is a compatible subset of $\Psi(Q)_{\geq -1}$ and $\alpha \in \cone(C)$. Furthermore, it is clear that there are only finitely many such compatible sets $C$ as $Q(\mathcal{E})$ has finitely many positive roots.
\end{proof}

\begin{remark}
Note that \cite[Theorem 8.1.7]{IOTW} can also be deduced from the proposition above, Example \ref{qes-ex}, and Theorem \ref{domains-exceptional-thm}.
\end{remark}

\end{document}